# Zur Gültigkeit der Annahmen des klassischen Sammelbilderproblems

Niklas Braband[1], Sonja Braband[1] und Malte Braband[2]

[1]Gymnasium Neue Oberschule, Braunschweig, Deutschland
`{niklas.braband|sonja.braband}@no-bs.de`
[2]Technische Universität Braunschweig, Deutschland
`m.braband@tu-braunschweig.de`

**Zusammenfassung.** Im klassischen Sammelbilderproblem wird angenommen, dass die Sammelbilder gleichverteilt und unabhängig sind, d. h. alle Bilder kommen gleich häufig vor und werden zufällig auf die Päckchen verteilt. Wir können sowohl statistisch als auch analytisch nachweisen, dass insbesondere die Annahme der Zufälligkeit bzw. Unabhängigkeit in der Praxis nicht erfüllt ist. Dazu haben wir den Herstellungsprozess der Sammelbilder systematisch mit kombinatorischen Methoden untersucht. Wir können auch Maße dafür angeben, wie groß die Abweichung vom Zufall ist.

Damit sind die Ergebnisse des klassischen Sammelbildermodells in der Praxis nicht mehr gültig, aber wir können zeigen, dass sie eine Abschätzung nach oben darstellen. Dies bedeutet insbesondere, dass kein Betrug durch den Hersteller vorliegt, sondern vielmehr eine Bevorteilung der Sammler.

**Abstract.** The solution of the classical Coupon Collector's Problem is based on the assumptions that all stickers are independently and uniformly distributed. We can prove statistically as well as analytically that in particular the assumption of independence is not fulfilled in the field. To achieve this result we have systematically modeled and analyzed the production process of the stickers by combinatorial methods. We have also derived measures for the deviation from independence.

Thus the results for the classical Coupon Collector's Problem are not valid for practical applications anymore, but we can show that they constitute an a upper bound. This means in particular that in practice the deviation is advantageous for the collectors, not a fraud of the vendors.

# 1    Einführung

Sammelbilder gibt es zu vielen Themen, z.B. Fussballbilder zur WM. Und es gibt viel Diskussionen über sie, z.B. viele Doppelte, manche Bilder bekommt man gar nicht, einige Sammler vermuten sogar Betrug der Hersteller. Wir haben begonnen, uns für das Thema zu interessieren [1], als Schweizer Wissenschaftler eine optimale Sammelstrategie [2] vorgeschlagen haben, die auch in der Presse breit veröffentlicht wurde.

Wikipedia [3] schreibt "Das Sammelbilderproblem .... befasst sich mit der Frage, wie viele zufällig ausgewählte Bilder einer Sammelbildserie zu kaufen sind, um eine komplette Bildserie zu erhalten". D.h. es geht darum, eine bestimmte Anzahl B von Bildern vollständig zu sammeln. Man kauft die Bilder normalerweise nicht einzeln, sondern in Päckchen, in denen je P Bilder drin sind. Bei Panini war z.B. in Deutschland zur WM 2014 B=640 und P=5.

Im klassischen Sammelbilderproblem werden die folgenden Annahmen getroffen:

A1.  Die Bilder sind rein zufällig auf die Päckchen verteilt. D.h. der Hersteller mischt bei der Herstellung ordentlich.
A2.  In einem Päckchen kommt kein Bild doppelt vor.
A3.  Alle Bilder kommen gleich häufig vor. D. h. der Hersteller betrügt nicht durch absichtliche Verknappung von Bildern.

Sardy und Velenik [2] haben behauptet, dass es optimal ist, zuerst ein Display zu kaufen (häufig auch Box genannt) zu kaufen (bestehend aus 100 Päckchen mit je 5 Bildern). Ein Display ist meistens billiger, als die Päckchen einzeln zu kaufen. Dabei gehen die Forscher davon aus, dass in einem Display alle Bilder verschieden sind. Guckt man aber im Internet bei Amazon die Kommentare zu den Displays an, so gibt es Beschwerden über viele doppelte Bilder. In den uns überlassenen Displays der Serie Match Attax von Topps gab es auch Doppelte. Auch in den Angeboten wird nicht damit geworben, dass alle Bilder in einem Display verschieden sind. Manche meinen, dass es trotzdem weniger Doppelte gäbe als beim Kauf der einzelnen Päckchen.

In der Literatur sind viele Ergebnisse über das klassische Sammelbilderproblem [2,3,4,5,6,7] zu finden, die alle auf denselben Annahmen A1-A3 oben basieren. A2 wird von allen Herstellern garantiert, der Effekt ist aber gering [3], so dass er meistens vernachlässigt wird. Nur vereinzelt gibt es Zweifel an den anderen Annahmen, wie z. B. in einer Analyse, die zur Fußball-WM über das Internet gemacht wurde [9]. Dort wurden Zweifel an der Gleichverteilung genährt, d. h. A3 in Frage gestellt. Aber letztendlich wäre diese Abweichung auch durch Fertigungsprobleme erklärbar und die Aussagekraft einer Datensammlung über das Internet ist auch fragwürdig.

Die Hersteller behaupten in Presseberichten [10], aber auch in ihren eigenen Medien [11], dass die Verpackung der Bilder in die Päckchen zufällig erfolgt. Wir haben uns Gedanken gemacht, wie der Hersteller es schafft, zwei sich eigentlich widersprechende Ziele unter einen Hut zu bekommen: Einerseits verspricht er, dass in einem Päckchen niemals Doppelte vorkommen (A2), andererseits, dass alle Bilder zufällig und gleichverteilt vorkommen (A1 bzw. A3). Dies entspricht Ziehen ohne Zurücklegen. Wie kann er das praktisch umsetzen? Wenn alle Bilder gut gemischt würden, könnte man nicht sicherstellen, dass in Päckchen keine Doppelten vorkommen oder man müsste das mühsam nachträglich kontrollieren. Aber zur Analyse müsste man entweder sehr viele Bilder haben oder den Herstellungsprozess kennen, den die Hersteller nach Presseberichten aber geheim halten.

## 2    Experimente

Als der Hersteller Topps im Herbst 2014 zusammen mit Penny die Bundesliga-Bilder wieder angeboten hat, haben wir mitgesammelt, aber eher um die bekannten Formeln zu bestätigen. Bei Penny bekam man kostenlos für je 10€ Einkaufswert ein Päckchen mit 5 Bildern gratis dazu. Insgesamt waren B=300 Bilder zu sammeln und man konnte K=50 Bilder nachkaufen. Nach den klassischen Formeln wären im Mittel etwa 535 Bilder zu kaufen, bevor man nachkaufen kann. Bis zum Februar 2015 hatten wir aber nur 34 Päckchen erhalten und damit 129 Bilder im Album gesammelt. Freundlicherweise hat uns die Firma Topps noch einmal 50 Päckchen zur Verfügung gesponsert. Nach Öffnen dieser 50 Päckchen haben uns dann nur noch 33 Bilder gefehlt. Schon nach 78 Päckchen oder 390 Bildern hätten wir nachkaufen können.

Wie man schon in Figur 1 erkennen kann, ist die Abweichung, die sich erst bei den von Topps nachgelieferten Bildern ergeben hat, sehr deutlich. Figur 2 zeigt die Anzahl der Doppelten pro Päckchen sowie einem gleitenden Durchschnitt (10 Päckchen) als Trendlinie. Wir hätten erwartet, dass es immer mehr Doppelte werden. Bei unserem Album scheint die Anzahl aber konstant zu bleiben.

Wir haben uns daher die Frage gestellt, ob wir einfach nur großes Glück gehabt haben oder ob es vielleicht eine andere plausible Erklärung gibt. Dazu rechnen wir aus, mit welcher Wahrscheinlichkeit ein solches Ergebnis vorkommt, wenn die Bilder zufällig verteilt wären.

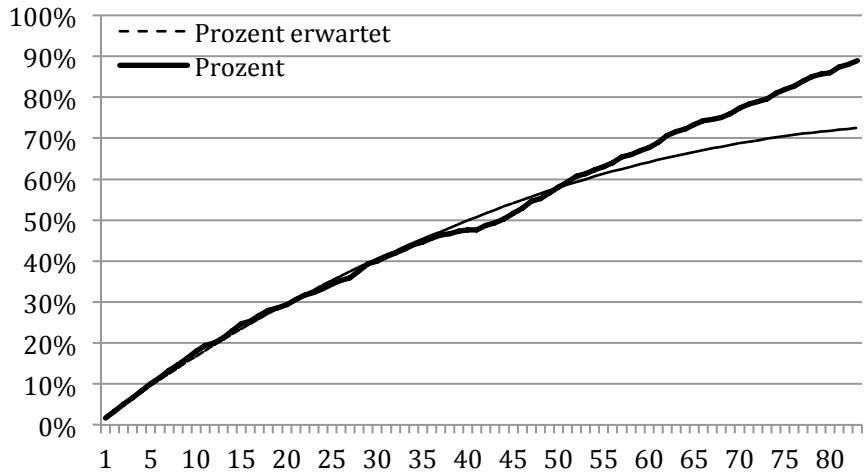

**Fig. 1.** Verlauf des Bundesliga-Sammelexperiments

Mit Hilfe der geometrischen Verteilung haben wir eine Standardabweichung von ca. 31 berechnet, bis man beim Bundesliga-Album nachkaufen kann. Dies bedeutet, dass unser Ergebnis fast 5σ vom Erwartungswert entfernt liegt. Wenn die Verteilung eine Normalverteilung wäre, so wäre die Wahrscheinlichkeit, dass eine solche Abweichung bei zufälliger Verteilung der Bilder auftritt, kleiner als 1:1,000,000.

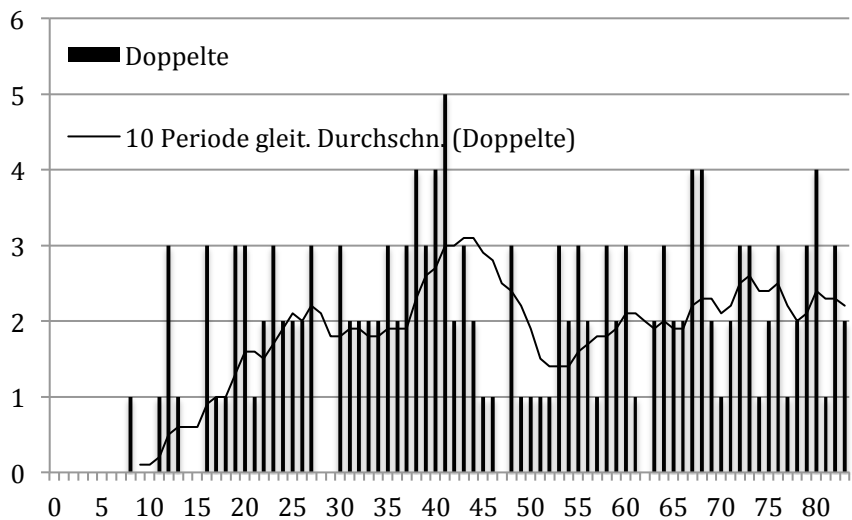

**Fig. 2.** Verlauf der Doppelten in den Päckchen

Da die exakte Verteilung schwer zu bestimmen ist, haben wir die Verteilung simuliert, und zwar zuerst mit 10 Millionen Versuchen. Dabei war das kleinste Ergebnis 405 Bilder, das auch nur einmal aufgetreten ist, siehe Figur 4. Aus diesen Ergebnissen konnten wir mit R die Schiefe der Verteilung schätzen und das Ergebnis 0,21 weist auf eine rechtsschiefe Verteilung hin. D. h., die Verteilung ist nicht ganz symmetrisch und fällt links steiler ab als rechts. Dies bedeutet aber, dass kleine Ergebnisse noch unwahrscheinlicher sind als bei der Normalverteilung. Auch bei einer Simulation mit insgesamt 50 Millionen ist ein derart kleines Ergebnis nicht aufgetreten. Insgesamt heißt das, dass das Ergebnis unseres Experiments mindestens so unwahrscheinlich war, wie im Lotto zu gewinnen.

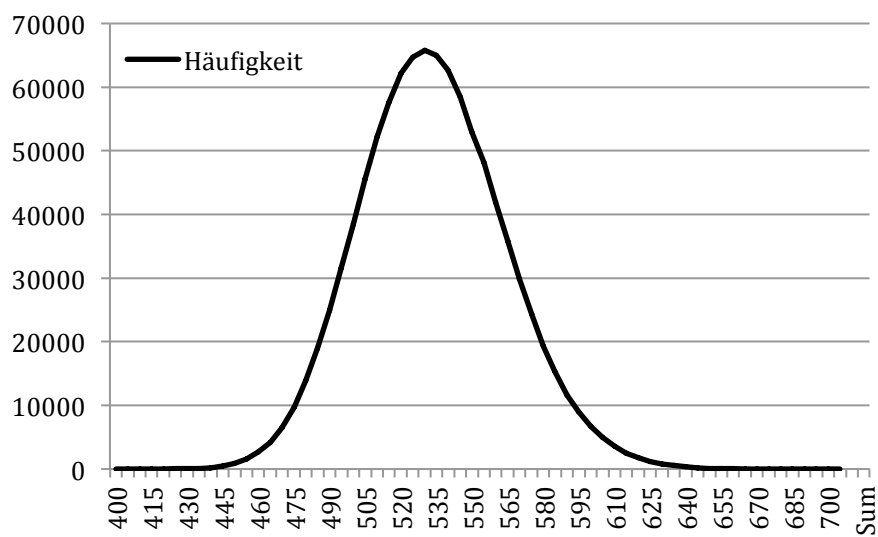

**Fig. 3.** Häufigkeitsverteilung der Bilder bei 10 Millionen Simulationen

Wir haben daraufhin den Hersteller Topps kontaktiert, der uns geantwortet hat [12]:"Dann hatten Sie bei den Päckchen, die wir Ihnen zugesendet haben, wohl ziemlich viel Glück. Wie schon in unserer Kommunikation zuvor beschrieben, werden die Sticker in den Päckchen absolut zufällig gemischt." Dies spornte uns an, uns noch etwas genauer mit dem Thema zu beschäftigen.

Aber wir standen vor dem Problem, dass ein Display mit 300-500 Bildern je nach Sammelserie ziemlich teuer ist und wir uns gefragt haben, wie wir an weiteres Datenmaterial kommen können. Dann erinnerten wir uns an hunderte von Bewertungen im Amazon Marketplace zu Käufen von Panini Displays zur Fußballweltmeisterschaft 2014. In vielen dieser Kritiken gab es auch Angaben zur Anzahl der Doppelten in den

Displays. Wir bezogen nur Angaben in die Auswertung ein, die konkrete Zahlenwerte genannt haben. Diese variierten von null Doppelten bis 330 Doppelten (bei 500 Bildern). Der Mittelwert lag bei ca. 130 Doppelten. Unter den klassischen Annahmen [3] hätten wir etwa 153 erwartet. Aber natürlich bleiben Zweifel an der Glaubwürdigkeit der Aussagen zu Doppelten aus den Amazon-Bewertungen.

Daher entschlossen wir uns, bei eBay doch noch zwei gelbe Original-Displays günstig kaufen. Wir zählten 73 bzw. 98 Doppelte pro Display. Später haben wir noch ein grünes Original-Display, welches für den Verkauf in Osteuropa bestimmt war, dazu gekauft. Hier traten 111 Doppelte im Display auf, aber die Bilder waren im Vergleich zu den gelben Displays schlechter gemischt, z.B. haben sich Folgen von Doppelten deutlich sichtbar wiederholt. Wir haben dann durch Simulation ermittelt, dass in 95% der Fälle zwischen 138 und 168 Doppelte im Display sein müssten, wenn die Bilder zufällig verteilt wären. Das bedeutet, dass sehr viele der bei Amazon angegebenen Doppelten und auch unsere Ergebnisse stark signifikant gegen die Hypothese der Zufälligkeit sprechen.

Schließlich kontaktierten wir Panini direkt. Glücklicherweise sponserten sie uns Displays mit je 50 Päckchen zu 6 Bildern der Serie Amici Cucciolotti, eines Tierbilderalbums. Dabei müssen insgesamt B=576 Bilder gesammelt werden, das heißt, das Album ist etwa gleich groß wie das WM-Album. Die Hälfte der Bilder sind Glitzer- oder Sonderbilder. Schon beim Auspacken der ersten Displays gab es Überraschungen, z. B. lagen bei den Päckchen lagen immer die Glitzer- und Sonderbilder obenauf und bei allen 10 Displays gab überhaupt keine Doppelten.

Unter den klassischen Annahmen [3] hätten wir im Mittel etwa 66 Doppelte pro Display erwartet. Die Wahrscheinlichkeit, bei zufälliger Mischung keine Doppelten zu erhalten, beträgt ungefähr $10^{-42}$. Mit einer Simulation haben wir ermittelt, dass mit 95% Wahrscheinlichkeit zwischen 54 und 78 Doppelte hätten auftreten müssen. Auch ist uns aufgefallen, dass es sehr lange Reihen aufeinander folgender Bilder gibt. Per Simulation haben wir ermittelt, dass mit 95% Wahrscheinlichkeit keine Folgen von mehr als 12 Bildern auftreten sollten, wenn die Bilder zufällig gemischt worden wären. Bei 8 von 10 Displays traten aber längere Folgen auf. Eine zufällige Mischung können wir mit diesen Ergebnissen daher so gut wie sicher ausschließen.

Dies bedeutet aber, dass beim Kauf von Displays die Annahme der zufälligen Verteilung verletzt ist (und zwar sowohl bei Panini als auch bei Topps), auf der alle Berechnungen zum klassischen Sammelbilderproblem beruhen! Wenn der Sammler die Päckchen in kleinen Mengen kauft, bemerkt er das nicht, so wie wir es bei den ersten 34 Päckchen unseres ersten Experiments (siehe Figur 1) nicht bemerkt haben. Dadurch, dass wir zu unterschiedlichen Zeitpunkten gekauft und Päckchen aus unter-schiedlichen Displays erhalten haben, sah das Ergebnis für uns zufällig aus. Man muss aber auch einräumen, dass hier kein Betrug, sondern eine Bevorteilung der Sammler durch die Hersteller vorliegt.

Eine Rückfrage bei Panini Deutschland [13] ergab, dass es bei der Serie Amici Cucciolotti gewollt ist, dass keine Doppelten im Display vorkommen. Diese Antwort und die verblüffenden Ergebnisse haben uns motiviert, die Herstellung der Displays genauer unter die Lupe zu nehmen.

## 3  Analyse des Herstellungsprozesses

Wir haben uns Gedanken gemacht, wie der Hersteller es schafft, zwei sich eigentlich widersprechende Ziele unter einen Hut zu bekommen: Einerseits verspricht er, dass in einem Päckchen niemals Doppelte vorkommen, andererseits, dass alle Bilder zufällig d. h. gleichverteilt vorkommen. Dies entspricht Ziehen ohne Zurücklegen. Wie kann er das praktisch umsetzen? Wenn alle Bilder gut gemischt würden, könnte man nicht sicherstellen, dass in Päckchen keine Doppelten vorkommen, oder man müsste das mühsam nachträglich kontrollieren.

Allerdings hatten wir zu Topps zunächst keine weiteren Informationen gefunden. Später wies uns Topps auf einen Artikel [10] hin, der allerdings weniger detailliert als die Informationen zu Panini war, bei denen es seit 2014 Filme über die Produktion gibt [11]. In Paninis Firmengeschichte kann man nachlesen, dass das Mischen anfangs tatsächlich mit viel Handarbeit verbunden war: „… es war notwendig, die Bilder sehr lange zu mischen, bevor sie eingetütet wurden, um die völlige Zufälligkeit zu garantieren; es durfte kein Bild doppelt in einem Tütchen sein … Die Sammelbilder wurden anfangs mit Hilfe eines Butterfasses gemischt, inspiriert durch die Maschine für die Butterherstellung und mit dem Rad für die Lottozahlen." Später wurde es notwendig, die Produktion zu automatisieren und Panini verweist auf ihre Fifimatic genannte, patentierte Verpackungsmaschine.

Wir stellen hier ein Modell des Herstellungsprozesses vor, der auf den Firmendarstellungen [11] sowie einer Patentschrift zur Fifimatic [14] beruht:

1. Alle B Bilder werden auf sogenannte Quadrotten gedruckt. Das sind große Bögen, die 4n Sammelbilder enthalten. Jedes Bild kommt genau einmal auf einer Quadrotte vor und es gibt genau B/4n verschiedene Quadrotten. Für das WM-Album ist n=5 und das Tieralbum n=6.

2. Die Quadrotten werden gemischt und übereinander gestapelt. Sie werden in Einzelbilder zerschnitten, aber die Struktur der Quadrotten wird beibehalten.

Die sogenannte Fifimatic ist eine automatische Verpackungsmaschine, die mit den Sammelbildern bestückt wird. Sie besteht im wesentlichen aus einer rechteckigen Anordnung von viermal n Magazinen für Einzelbilder-Stapel (nach [11] variiert n in der Praxis zwischen 4 und 7, im Patent [14] ursprünglich nur 4). In Figur 4 sieht man die Seitenansicht der Fifimatic aus [14], hier bezeichnet Nr. 1 und 2 die Seiten der Magazine. Die Bilderstapel (Nr. 3 in Figur 4) werden so eingeordnet, dass die Struktur der Quadrotten erhalten bleibt. D. h. wir können unser Modell auch so beschreiben, als ob die Fifimatic direkt mit den Quadrotten arbeitet. Unter den Magazinen laufen 4 Fließbänder (Nr. 12 in Figur 4).

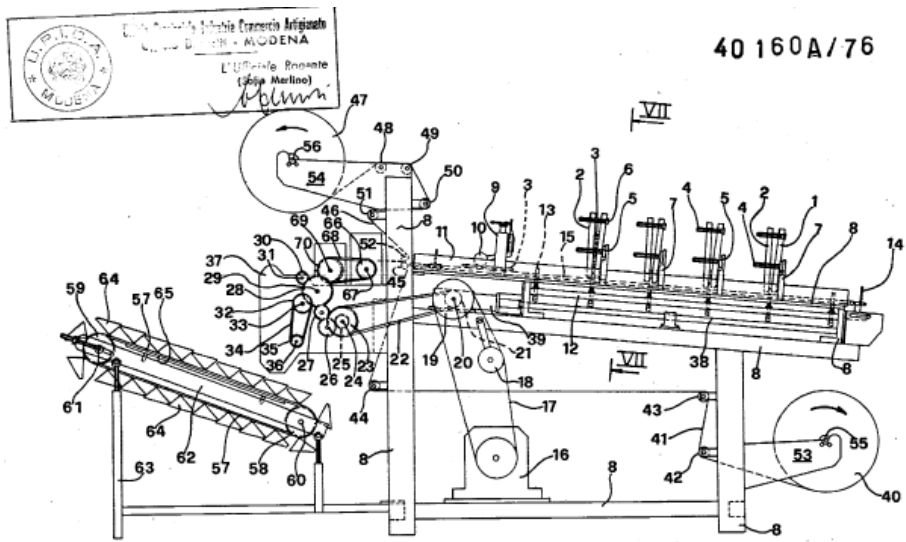

**Fig. 4.** Seitenansicht der Fifimatic nach [14]

Die eigentliche Produktion der Sammelbilderpäckchen verläuft nun weiter wie folgt:

3. In jedem Takt wird von jedem Magazin ein Bild nach unten auf die Fließbänder abgelegt, d. h. insgesamt wird in jedem Takt die unterste Quadrotte auf die Fließbänder gelegt.

4. Dann wird das Fließband eine Position weiter gefahren und die nächste Quadrotte abgelegt.

Am Ende entstehen 4 parallele Stapel mit je n Karten, die in Päckchen verpackt werden. Etwas mathematischer ausgedrückt, gibt es q=B/4n unterschiedliche Quadrotten $Q_1, Q_2, ... Q_q$. Jedes Sammelbild x aus {1, ..., B} hat eine eindeutige Position auf genau einer Quadrotte, d. h. jedes x entspricht genau einer Position $Q_k(i,j)$, mit i=0, ...n-1 und j=0,...,3 (die Spalten der Quadrotte entsprechen den Fließbändern). Wenn für

$y=Q_r(s,t)$ gilt, dann gilt $x=y$ genau dann wenn $r=k$ und $s=i$ und $t=j$ gilt. Zwei verschiedene Sammelbilder x und y können nur in ein Päckchen gelangen, wenn ihren Quadrotten in derselben Spalte stehen, d. h. wenn $j=t$ gilt.

Das Ergebnis des Mischens der Quadrotten kann man als eine Folge von Quadrotten $R_1, R_2,...$ auffassen, wobei jedes $R_i$ zufällig aus der Menge aller Quadrotten $\{Q_1, Q_2,... Q_q\}$ gewählt wird. Jetzt betrachten wir ein beliebiges Bild $R_k(i,j)$ aus der Produktion und untersuchen, mit welchen Bildern es in ein Päckchen verpackt werden kann. Aufgrund der Konstruktion der Fifimatic kann das Bild nur mit den Bildern $R_{k-n}(i+n,j)$, $R_{k-(n-1)}(i+n-1,j)$, ... $R_{k-1}(i+1,j)$, $R_{k+1}(i-1,j)$, ... $R_{k+n-1}(i-(n-1),j)$ oder $R_{k+n}(i-n,j)$ in ein Päckchen geraten (wobei die Indizes modulo n zu verstehen sind). Da jedes Bild eine eindeutige Position auf der Quadrotte hat, müsste erstens zweimal dieselbe Quadrotte unter den $R_{k-n}, ... R_{k+n}$ gezogen werden und zweitens müsste der Zeilenindex i des betrachteten Bildes gleich einem der Zeilenindizes der anderen Bilder sein, also gleich i+n, ..., i+1, i-1,..., i-n. Dies ist aber selbst dann unmöglich, wenn alle $R_{k-n}, ... R_{k+n}$ identisch wären, denn jedes Bild hat einen eindeutigen Zeilenindex, hier i, und der ist ungleich allen anderen Indizes i+n, ..., i+1, i-1,..., i-n in der Folge. Damit haben wir bewiesen:

**Satz 1:** Wenn jedes Bild eine eindeutige Position auf einer Quadrotte hat, d. h. $x=Q_k(i,j)$, und jedes Bild in die zugehörige i-te Zeile der Fifimatic einsortiert wird, dann gibt es keine doppelten Bilder in mit der Fifimatic produzierten Päckchen.

Figur 5 zeigt das Prinzip anhand des WM-Albums. Hier gibt es 32 verschiedene 4x5-Quadrottten. Diese werden gemischt und entsprechend in die Fifimatic einsortiert. Durch das Weiterschieben der Fließbänder um jeweils eine Position gelangen immer ein Bild aus einer ersten, zweiten, dritten, vierten und fünften Zeile einer Quadrotte in ein Päckchen. Darunter können dann keine doppelten sein, wenn die Bilder immer in die entsprechende Zeile der Fifimatic einsortiert werden, selbst wenn alle Quadrotten identisch wären, da jedes Bild eine eindeutige Position auf einer Quadrotte besitzt.

Die systematische Vermeidung von Doppelten hat allerdings ihren Preis: bei zufälliger Mischung der Sammelbilder gibt es $\binom{B}{P}$ verschiedene Päckchen. Wenn man die Fifimatic aber genauso befüllt, wie die Bilder auf den Quadrotten angeordnet sind, so können keine Bilder aus unterschiedlichen Spalten in ein Päckchen gelangen und in jedem Magazin nur B/4n unterschiedliche Bilder, d. h. bei zufälliger Mischung der Quadrotten gäbe es nur noch $4\left(\frac{B}{4n}\right)^n$ verschiedene Päckchen. Beim WM-Album mit B=640, P=5 und n=5 gibt es bei zufälliger Mischung etwa $881\times10^9$ verschiedene Päckchen, aber aufgrund der Fifimatic-Systematik kommen davon höchstens etwa

$134 \times 10^6$ vor, und dies auch nur unter Annahme der zufälligen Mischung der Quadrotten.

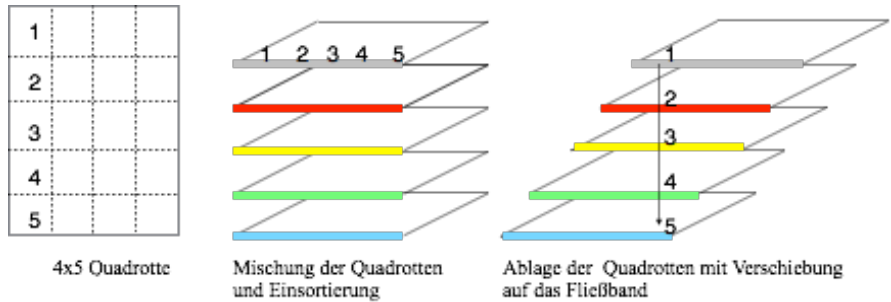

**Fig. 5.** Modell des Herstellungsprozesses am Beispiel des WM-Albums

Wir können damit behaupten

**Satz 2:** Beim Herstellungsprozess der Sammelbilder mit dem Fifimatic-Prinzip werden die Bilder nicht zufällig auf die Päckchen verteilt.

Die Abweichung betrifft vor allem die Annahme der Unabhängigkeit: Ist bekannt, dass Bild x=Qk(i,j) in einem Päckchen vorkommt, so kann kein anderes Bild y=Qr(s,t) in dem Päckchen vorkommen mit i=s oder j≠t. Dies bedeutet, dass die Wahrscheinlichkeit für das Auftreten des nächsten Bildes y abhängig ist vom aktuellen Bild x. Wenn die Quadrotten zufällig gemischt wären, könnten genau q=B/4n Bilder mit einer derselben Wahrscheinlichkeit auftreten, für alle anderen B-q Bilder ware die Auftretenswahrscheinlichkeit Null.

## 4 Doppelte in Displays

Jetzt bleibt noch die Frage zu klären, wie man es im Herstellungsprozess schaffen kann, dass keine oder weniger Doppelte in Displays vorkommen als eigentlich erwartet. Bei den von uns betrachteten Sammelserien gab es entweder sehr große Displays mit weniger Doppelten als erwartet (D=500 Bilder im Display im Vergleich zu B=640 Bildern im WM-Album, D=250 Bilder im Vergleich zu B=300 Bildern im Bundesliga-Album) oder es gab sogar Displays ohne Doppelte (D=300 Bilder im Vergleich zu B=576 Bildern im Tieralbum). Unsere Vermutung ist, dass dies mit der Mischung der Quadrotten zusammenhängt. Im Firmenfilm [11] kann man sehen, dass Quadrotten aus q artreinen Stapeln von beiden Seiten auf ein Fließband fallen. Dies scheint zu bestätigen, da dort eher wenig gemischt wird.

Daher betrachten wir jetzt Folgen $R_1$, $R_2$,... von Quadrotten und nehmen etwas vergröbert an, dass immer eine bestimmte Anzahl von Quadrotten in ein Display verpackt wird. Wir nehmen zuerst an, dass die Quadrotten gar nicht gemischt warden, sondern sich immer wieder die Folge $Q_1$, $Q_2$,... $Q_q$ wiederholt. Ist dann das Display kleiner als die Gesamtzahl der Sammelbilder (D<B), so gelangt immer eine gewisse Anzahl von Quadrotten in ein Display und es kommen in den Displays keine Doppelten vor. Auch wenn man jetzt in der Folge etwas mischt, z. B. an ein paar Plätzen Quadrotten vertauscht, so wird es bei kleinen Displays immer noch keine Doppelten geben und auch bei größeren Displays nur mit einer geringen Wahrscheinlichkeit. Der Firmenfilm legt die Vermutung nahe, dass durch die spezielle Konstruktion der Mischmaschine nur wenig und dann nur innerhalb der Folge $Q_1$, $Q_2$,... $Q_q$ gemischt wird. Wenn man in dieser Folge stärker oder wirklich zufällig mischt, so ergeben sich mehr Doppelte in den Displays. Im ungünstigsten Fall kann man erwarten, dass die Hälfte der Bilder in einem Display doppelt vorkommt. Dies passt sehr gut zu den Bewertungen im Amazon Marketplace, wo nur in ganz wenigen Fällen mehr als die Hälfte der Bilder als doppelt berichtet wurden. Diese Berichte könnten aber fehlerhaft sein. Unser Modell kann aber zumindest erklären, wie der Hersteller die Doppelten in den Displays steuern oder garantieren kann, dass keine Doppelten vorkommen.

Um zu überprüfen, ob unser Modell die Wirklichkeit angemessen beschreibt, haben wir zehn Displays aus der Tierbilder-Serie ausgewertet. Jedes Display (300 Bilder) selbst hatte keine Doppelten. Wir haben uns jeweils zwei Displays angesehen (zusammen 600 Bilder) und die Doppelten manuell ausgezählt. Bei B=576 können also zwischen 24 und 300 Doppelte vorkommen. Mit einer Simulation haben wir ermittelt, dass mit 95% Wahrscheinlichkeit zwischen 144 und 168 Doppelte vorkommen müssten, wenn die Bilder wirklich zufällig in den Displays gemischt wären. Würde unser Modell zutreffen, so müsste eine andere Verteilung zu beobachten sein, da auch sehr wenige oder sehr viele Doppelte vorkommen müssten.

|     | 217 | 218 | 219 | 220 | 221 | 526 | 530 | 531 | 533 |
|-----|-----|-----|-----|-----|-----|-----|-----|-----|-----|
| 216 | 136 | 166 | 140 | 127 | 24  | 148 | 96  | 141 | 143 |
| 217 |     | 146 | 173 | 151 | 156 | 135 | 149 | 121 | 168 |
| 218 |     |     | 169 | 148 | 141 | 144 | 136 | 158 | 175 |
| 219 |     |     |     | 152 | 170 | 144 | 155 | 171 | 162 |
| 220 |     |     |     |     | 157 | 130 | 149 | 193 | 210 |
| 221 |     |     |     |     |     | 156 | 210 | 159 | 155 |
| 526 |     |     |     |     |     |     | 143 | 24  | 90  |
| 530 |     |     |     |     |     |     |     | 170 | 122 |
| 531 |     |     |     |     |     |     |     |     | 229 |

**Tab. 1.** Doppelte in zwei Displays der Tierbilder-Serie

Die Ergebnisse in Tabelle 1 sprechen klar für unser Modell, denn es kommt eine breite Verteilung von Doppelten vor. Bei 45 Paarvergleichen hätten wir nur wenige Ausreißer erwartet, aber im Versuch gab es 27 signifikante Abweichungen, von der Minimalzahl von Doppelten (24) bis zu 229. Die Zahlen in den Tabellenköpfen geben die drei letzten Ziffern der Seriennummer des Displays an. Auf jedem Display war zusätzlich noch eine Nummer von 11 bis 15 aufgeklebt. Wir vermuten, dass sie die

Nummer der Fifimatic angibt, von der die Bilder stammen. Hier ergeben sich bei allen Paaren mit gleicher Nummer signifikante Abweichungen.

Dies bedeutet insgesamt, dass unser Modell den Herstellungsprozess stochastisch gut erklären kann, insbesondere warum es zu weniger Doppelten in den Displays kommt. Dieser Effekt führt insgesamt zu einer Bevorteilung der Sammler, es liegt kein Betrug seitens der Hersteller vor. Dies bedeutet aber auch, dass die Ergebnisse, die unter den Annahmen des klassischen Sammelbilderproblems hergeleitet wurden, in diesem Fall nicht mehr gültig sind, sondern nur noch Abschätzungen nach oben darstellen.

## 5 Zusammenfassung

Wir haben in dieser Arbeit sowohl statistisch als auch analytisch nachweisen, dass insbesondere die klassische Annahme A1 der Zufälligkeit bzw. Unabhängigkeit im Sammelbilderproblem in der Praxis nicht erfüllt ist. Dazu haben wir den Herstellungsprozess der Sammelbilder systematisch mit kombinatorischen Methoden untersucht und ein Modell für den Herstellungsprozess aufgestellt.

Wir haben große Abweichungen vom Zufall gefunden, z. B. kommen beim Panini WM-Album mindestens zehntausendmal weniger Möglichkeiten für Päckchen vor als bei zufälliger Mischung. Diese Abweichungen wurden vom Hersteller auch teilweise eingeräumt, z. B. ist bei manchen Serien gewollt, dass in den Displays keine Doppelten vorkommen.

Damit sind die Ergebnisse des klassischen Sammelbildermodells in der Praxis nicht mehr gültig, denn es kommen beim Kauf von Displays wesentlich weniger Doppelte vor als von der Theorie vorhergesagt. D. h. die klassischen Ergebnisse stellen nur noch obere Schranken dar. In einer anderen Arbeit [15] haben wir erste Ansätze aufgezeigt, wie man den Effekt von Displays in den klassischen Formeln berücksichtigen könnte.

Ein Sammler, der seine Päckchen mehr oder weniger einzeln am Kiosk kauft, wird diesen Effekt allerdings nicht bemerken. Dies bedeutet, dass kein Betrug durch den Hersteller vorliegt, wie häufig in den Medien behauptet, sondern sogar eine systematische Bevorteilung der Sammler durch den Hersteller. Dieser Effekt ist allerdings nur bei Kauf eines einzigen Displays wirksam, beim Kauf eines zweiten Displays besteht ein hohes Risiko, viele Doppelte zu bekommen.

## Referenzen